\documentclass[oneside,a4paper]{article}

\date{as of: 1 July 2002}% see end of file for history

\title{Algebras of generalized functions\\
  through sequence spaces algebras.\\
  Functoriality and associations}

\author{Antoine Delcroix\\\small\tt Antoine.Delcroix@univ-ag.fr
\and Maximilian F. Hasler%\thanks{
\\\small\tt MHasler@martinique.univ-ag.fr
\and Stevan Pilipovi\'c\\\small\tt pilipovic@im.ns.ac.yu
\and Vincent Valmorin\\\small\tt Vincent.Valmorin@univ-ag.fr}

\newcommand\comment[1]{}% comments are not shown if this command is defined
\usepackage{amsmath,amsfonts,amssymb,bm}% bm for bold greek letters

\makeatletter
\let\over\@@over% bug in amsmath

\catcode176=13 \def^^b0{\degres}%\ensuremath{{\,^\circ}}
\catcode178=13 \def^^b2{\ensuremath{^2}}
\catcode183=13 \def^^b7{\ensuremath{\cdot}}

\newcommand\A{{\mathcal A}}
\newcommand\B{{\mathcal B}}
\newcommand\C{{\mathbb C}}
\newcommand\CC{{\mathcal C}}
\newcommand\CEP{$\left(\mathcal C,\mathcal E,\mathcal P\right)$}
\ifx\undefined\comment
\newcommand\comment[1]{\mbox{}\marginpar{\flushleft\tiny#1}}
\fi

\newcommand\D{{\mathcal D}}
\newcommand\DD{{\text{\textbf{\textsl D}}}}
\newcommand\e{\varepsilon}
\newcommand\E{{\mathcal E}}

\newcommand\F{{\mathcal F}}

\newcommand\G{{\mathcal G}}
\def\H{{\mathcal H}}% \H = hungarian long umlaut accent
\newcommand\I{{\mathcal I}}
\def\INT#1d{\int#1\,\mathrm d}
\newcommand\ind{\operatorname*{ind\,lim}}
\newcommand\impl{\mathop{~\Longrightarrow~}}
\newcommand\J{{\mathcal J}}
\newcommand\K{{\mathcal K}}
\newcommand\KK{{\mathbb K}}
\def\L{{\mathcal L}}%% ! was: L-bar !
\def\l{\lambda}%% ! was: l-bar !

\newcommand\lr[3]{\left#1#3\right#2}
\newcommand\Mid{~\Big|~}
\newcommand\N{{\mathbb N}}%I\!\!N}
\newcommand\nb[1]{\par\noindent\textbf{#1} }
\newcommand\ola{\overleftarrow}
\newcommand\olra{\overleftrightarrow}
\newcommand\ora{\overrightarrow}
\def\P{{\mathcal P}}
\newcommand\p[1]{\left(#1\right)}

\newcommand\R{{\mathbb R}}%I\!\!R}
\newcommand\set[1]{\left\{\,#1\,\right\}}

\newcommand\ultra[1]{{|\!|\!|\,#1\,|\!|\!|}} % |||\([^|]*\)|||

\newcommand\vp{\varphi}%  \v is check accent !!
\newcommand\veps{\varepsilon}
\newcommand\Z{{\mathbb Z}}

\newtheorem{theorem}{Theorem}

\newtheorem{definition}[theorem]{Definition}
\newtheorem{example}[theorem]{Example}
\newtheorem{lemma}[theorem]{Lemma}
\newtheorem{proposition}[theorem]{Proposition}
\newtheorem{propdef}[theorem]{Proposition--Definition}
\newtheorem{remark}[theorem]{Remark}
\newenvironment{proof}[1][Proof]{\nb{#1.}}{\hfill$\Box$\par}

\let\mydate\@date%\maketitle deletes \@date

\advance\footskip1cm

\numberwithin{equation}{section}
\begin{document}
\maketitle

\centerline{\bf Abstract}

Starting from a description of various generalized function algebras based on
sequence spaces, we develop the general framework for considering linear
problems with singular coefficients or non linear problems.  Therefore, we
prove functorial properties of those algebras and show how weak equalities, in
the sense of various associations, can be described in this setting.

\nb{AMS classification:} % http://www.ams.org/msc/46Axx.html
46A45 (Sequence spaces), % (including Köthe sequence spaces) [See also 46B45]
46F30 (Generalized functions for nonlinear analysis).
        % (Rosinger, Colombeau, nonstandard, etc.)
secondary:
46E10, % Topological linear spaces of continuous,
        % differentiable or analytic functions
%46A03, General theory of locally convex spaces
46A13, % Spaces defined by inductive or projective limits (LB, LF, etc.)
46A50, % Compactness in topological linear spaces; angelic spaces, etc.
% 46H05 General theory of topological algebras
46E35, % Sobolev spaces and other spaces of ``smooth'' functions,
% embedding theorems, trace theorems
% FOR SECOND PART :
% 46A22 extension and lifting of functionals and operators [See also 46M10]
46F05.% Topological linear spaces of test functions,
% distributions and ultradistributions [See also 46E10, 46E35]
%46F10 Operations with distributions
%46F15 Hyperfunctions, analytic functionals
%46F20 Distributions and ultradistributions as boundary values
%      of analytic functions [See also 30D40, 30E25, 32A40]

\section*{Introduction}

Schwartz~\cite{schw} proved that differential algebras of generalized
functions
% containing (embedded)  the delta distribution
with the ordinary product of continuous functions and containing the delta
distribution, do not exist.  But with the ordinary multiplication of smooth
functions, such algebras exist as it was proved by Colombeau. Nowadays the
theory of such algebras is well-established and it is affirmed through many
applications especially in nonlinear problems with strong singularities. We
refer to the books \cite{bia,nepisca,gros,ob} and to the numerous papers given
in the references.

We have shown in \cite{DHPV1} that Colombeau type algebras can be reconsidered
as a class of sequence space algebras and have given a purely topological
description of Colombeau type algebras.
%
%We will show that such algebras fit very well in the general
%theory of the well known sequence spaces forming appropriate algebras.
%
In fact, all these classes of algebras are simply determined by a locally
convex algebra $E$ and a sequence $r:\N\to\R_+$ (or sequence of sequences)
which serves to construct an ultrametric on subspaces of
% the sequence space
$E^\N$.  Such sequences are called weight sequences.

Distribution, ultradistribution and hyperfunction type spaces can be embedded
into corresponding algebras of sequences of this class.  This is done in
\cite{DHPV2}.
Note, the embeddings of Schwartz' spaces into the Colombeau algebra $\G$ are
very well known, but for ultradistribution and hyperfunction type spaces this
is quite different problem, especially because of multiplication of regular
enough functions (smooth, ultradifferentiable or quasianalytic), embedded into
corresponding algebras.

In this paper we continue to develop the foundations of our approach as well
as the general framework for considering various linear problems with singular
coefficients and nonlinear ones.

%We recall in Section  1 our construction of
%exponentially weighted sequence spaces forming algebras
%which correspond to an appropriate decreasing sequence.
%
We recall in Section~\ref{sect:basic} our construction of
algebras of sequence spaces, defined by a decreasing null sequence,
used as exponential weight.
Taking $r_n=1 / \log n$, $n\in\N$, we get
%can describe
the simplified and the full Colombeau algebras.

In Section~\ref{sect:scales} we introduce sequences of scales, construct new
algebras and relate them to known algebras as Egorov algebras~\cite{eg}
%, \CEP--algebras~\cite{CEP}
and asymptotic  algebras~\cite{asymp}.

This justifies to turn in Section~\ref{sect:functor} to nowadays classical
questions like functorial aspects of Colombeau type algebras~\cite{asymp,sca}
in order to apply the following scheme in standard applications: if a
classical differential problem for regular data has a unique solution such
that the map associating the solution to the initial data verifies convenient
growth conditions (with respect to the chosen scale of weights), then this
problem can be transferred to corresponding sequence spaces, where it also
allows for a unique solution. That way, differential problems with singular
data can be solved in such spaces {\em ad hoc}.

%Finally, it occurs frequently that exact solutions are not required; sometimes
%they do not exist at all.
Finally, exact solutions may not exist at all or, even more frequently,
may not be needed.
For this reason, in spaces of generalized functions
the notion of weak solutions has often be used, in the sense of different
types of associations.  These concepts can nicely be described in our
sequential approach, which is done in Section~\ref{sect:association}.  Indeed,
we give a generalized and unified scheme of a large number of tools of this
kind, including those which can already be found in various places in existing
literature~\cite[and others]{biacol1,ma1,CEP,nepisca}.

\section{The basic construction \protect\cite{DHPV1}}\label{sect:basic}

\subsection{Locally convex vector spaces and algebras}

Consider an algebra $E$ which is a locally convex vector space over $\C$,
equipped with an arbitrary set of seminorms $p\in\P$ determining its locally
convex structure. Assume that
\[
  \forall p\in\P~\exists \bar{p}\in\P~\exists C\in\R_+:
  \forall x,y\in E: p(x\,y) \leq C\,\bar{p}(x)\,\bar{p}(y) ~.
\]

Let $ r \in \R_+^\N $ be a sequence decreasing to zero. Put, for $ f\in E^\N $,
\[
  \ultra{f}_{p,r} := %\underset{n\to\infty}
  \limsup_{n\to\infty} p(f_n)^{r_n}~.
\]
This is well defined for any $f\in E^\N$, with values in
$
  \overline{\R }_+ := \R_+ \cup \{ \infty \}
$.
With this definition, let
\[
  \F_{\P,r} = \set{ f\in E^\N  \mid
  \forall p\in\P  : \ultra{f}_{p,r}< \infty } ~,
\]
%and
\[
  \K_{\P,r}=\set{ f\in E^\N \mid
  \forall p\in\P  : \ultra{f}_{p,r} = 0 } ~.~
\]
\begin{remark}
The space $E$ is given as domain of the elements $p\in\P$. When we
write $\F_{|\cdot|,r}$ in the sequence, it means $E=\C$ and
$\P=\set{|\cdot|}$.
\end{remark}

%Then the following holds:

\begin{proposition}\label{prop2}

\begin{enumerate}
\item $\F_{\P,r}$ is a (sub-)algebra of $E^\N $,
and $\K_{\P,r}$ is an ideal of $\F_{\P,r}$, thus
$\G_{\P,r}=\F_{\P,r}/\K_{\P,r}$ is an algebra.

\item  For every $p\in\P$, $%\vskip-8ex
\begin{array}[t]{rcl}
        d_{p,r}:~ E^\N \times E^\N &\to& \overline{\R }_+
~,\\                     (f,g)   &\mapsto& \ultra{f-g}_{p,r}
\end{array}$
\\
is an ultrapseudometric on $\F_{\P,r}$, and the family 
$\p{d_{p,r}}_{p\in\P}$ makes $\F_{\P,r}$ a topological algebra
(over $(\F_{|·|,r},d_{|·|,r})$).

\item For every $p\in\P$, $%\vskip-8ex
\begin{array}[t]{rcl}
        \widetilde{d}_{p,r} :~
        \G_{\P,r}\times\G_{\P,r} & \to& \R_+
~,\\
        ([f],[g]) & \mapsto& d_{p,r}(f,g)
\end{array}$\\
is an ultrametric on $\G_{\P,r}$, where $[f],[g]$ are the classes
of $f,g\in\F_{\P,r}$.\\[0.5ex]
The family of ultrametrics $\set{\smash{\widetilde{d}_{p,r}}}_{p\in\P}$
defines a topology, identical to the quotient topology, for which
$\G_{\P,r}=\F_{\P,r}/\K_{\P,r}$ is a topological algebra
over $\C_r=\G_{|·|,r}$.
\end{enumerate}
\end{proposition}

\begin{example}[Colombeau generalized numbers]\label{ex-numbers}~
The setting con\-si\-de\-red here is used to define rings of generalized
numbers.  For this, $E$ is the underlying field $\R$ or $\C$, and
$p=|\cdot|$ the absolute value. The resulting factor algebra
$\G_{|·|,r}$, with topology given by $\ultra·_{|·|,r}$, will be noted
$\R_r$ or $\C_r$.
As already explained in the introduction, for $r=1/\log$, we get the
ring of Colombeau's numbers $\overline\C$. More precisely, let
\begin{eqnarray}
  \forall n\in\N +2 : r_n=\frac{1}{\log n}~.
\label{*}
\end{eqnarray}
This gives back Colombeau's algebras of elements with polynomial growth modulo
elements of more than polynomial decrease, because%
\begin{eqnarray*}
  \limsup|x_n|^{1/\log n} < \infty
&\iff&  \exists C:\lim\sup|x_n|^{1/\log n}=C
\\
&\iff&  \exists B,\exists n_0, \forall n>n_0 :
   |x_n| \leq B^{\log n} = n^{\log B}
\\
&\iff&  \exists\gamma: |x_n| = o(n^\gamma) ~.
\end{eqnarray*}
%It is clear that there is actually equivalence
% between the first and last line.
On the other hand, $\limsup=0$ (for the ideal) corresponds to
% taking
$C=0$ and thus $\forall B>0$ and $\forall\gamma$ in the last lines.
\end{example}

\begin{example}
Take $ E=\CC^\infty(\Omega) $, $\P = \set{ p_\nu }_{\nu\in\N} $, with
$$
  p_\nu (f) := \sup_{|\alpha|\leq \nu,\,|x|\leq\nu}
  |f^{(\alpha)}(x)|
~,$$
and $r=\frac{1}{\log}$. Then, $\G_{\P ,r} = \F_{\P ,r} /\K_{\P ,r} $
with
\begin{eqnarray*}
  \F_{\P ,r} &=&
    \left\{ {(f_n)}_n \in {\CC^\infty(\Omega)}^\N
  ~\Big|~ \forall \nu\in\N:
  \limsup_{n\to \infty} p_\nu(f_n)^{1/\log n}< \infty \right \} ~,
\\
  \K_{\P ,r} &=&
  \left\{ {(f_n)}_n \in {\CC^\infty (\Omega)}^\N
  ~\Big|~ \forall \nu \in \N:
  \limsup_{n\to\infty} p_\nu(f_n)^{1/\log n}=0 \right\} ~.
\end{eqnarray*}
is just the simplified Colombeau algebra.
\end{example}

\begin{example}\label{Ex-Full-Colombeau}
With a slight generalization, we can consider the ``full'' Colombeau
algebra based on these definitions. Following Colombeau, let
\[
  \forall q\in \N :~~
  \A_{q} = \set{
  \phi\in\D( \R ^{s})
~\Big|~
  \int\phi=1
~~\text{and}~~
  \forall n\in\left\{ 1,\dots,q \right\} :\int x^n \phi=0
  }
\]
and $p_\nu$, $\nu\in\N$ as above. Then, for fixed $ \nu,N\in\N $ and
$\phi\in\A_N$ let
\begin{eqnarray*}
  \F_{\nu,N,\phi} &=&\set{ \p{f_{\phi_n}}_n\in E^\N \mid
        \ultra{\p{f_{\phi_n}}_n}_{p_{\nu},r}< N } ~,
\\
  \K_{\nu,N,\phi} &=&\set{ \p{f_{\phi_n}}_n\in E^\N \mid
        \ultra{\p{f_{\phi_n}}_n}_{p_{\nu},r}=0 } ~,
\end{eqnarray*}
where
$
  \phi_n = n^s \, \phi( n\,\cdot)
$.\\
(Here ${(f_{\phi_n})}_n$ are the ``extracted sequences'' of the elements
$ {(f_\phi)}_\phi \in E^{\D(\R)} $).\\
Now define
\begin{eqnarray*}
  \F  &=& \bigcap_{\nu\in\N }\F_{\nu} ~,~~
  \F_{\nu}= \bigcup_{N\in\N }\F_{\nu,N} ~,~~
  \F_{\nu,N}= \bigcap_{\phi\in\A_N}\F_{\nu,N,\phi} ~,
\\
  \K  &=& \bigcap_{\nu\in\N }\K_{\nu} ~,~~
  \K_{\nu}= \bigcup_{N\in\N }\K_{\nu,N} ~,~~
  \K_{\nu,N}= \bigcap_{\phi\in\A_N}\K_{\nu,N,\phi} ~.
\end{eqnarray*}
Then again, $\F$ is an algebra and $\K$ an ideal of $\F$, and
$\G=\F\,/\,\K$ is the ``full'' Colombeau algebra.
\end{example}

\subsection{Projective and inductive limits}
%*******************************************

We consider again a positive sequence $ r={(r_n)}_n\in(\R_+)^\N $ decreasing
to zero and we use the notations introduced above.

%If $p$ is a seminorm on a vector space $E$, we define for
%$f={(f_n)}_n\in E^\N$
%$$
%   \ultra f_{p,r} = \limsup_{n\to\infty} \big( p(f_n) \big)^{r_n}
%$$
%%with values in $\overline\R_+ = \R_+\cup\set\infty$.
%Denote $\tilde E^\N=\{f\in E^\N \mid \ultra f_{p,r} <\infty\}.$
%\\

Let $\left( E_\nu^\mu,p_{\nu\,}^\mu\right)_{\mu,\nu\in\N}$ be
a family of semi-normed algebras over $\R$ or $\C$ such that
\begin{align*}
        \forall\mu,\nu\in\N :~ E_\nu^{\mu+1} \hookrightarrow E_\nu^\mu
\text{ ~ and ~ }
        E_{\nu+1}^\mu \hookrightarrow E_\nu^\mu
\text{~ (resp. }
        E_\nu^\mu \hookrightarrow E_{\nu+1}^\mu ~)~,
\end{align*}
where $\hookrightarrow$ means continuously embedded.
%~(For the $\nu$ index we consider inclusions in the two directions.)~
%
%This implies that there exist constants $ C_\nu^\mu\in\R _+ $
%such that
%$ \forall\mu,\nu\in\N : p_\nu^\mu\leq C_\nu^\mu\,p_{\nu+1}^\mu $
%(on domains where both seminorms are defined),
%but without loss of generality one can take these constants equal to 1.
%The same holds for the $_\nu$ indices.
%
Then let
$\displaystyle
        \ola{E}
        =\projlim_{\mu\to\infty} \projlim_{\nu\to\infty} E_\nu^\mu
        =\projlim_{\nu\to\infty} E_\nu^\nu
$, (resp.
$\displaystyle
        \ora{E}
        =\projlim_{\mu\to\infty} \ind_{\nu\to\infty}E_\nu^\mu)
$.
Such projective and inductive limits are usually considered with norms
instead of seminorms, and with the additional assumption that in the
projective case sequences are reduced, while in the inductive case for
every $\mu\in\N$ the inductive limit is regular, i.e. a set
$A\subset{\ind\limits_{\nu\to\infty}}E_\nu^\mu$ is bounded iff
it is contained in some $E^\mu_\nu$ and bounded there.

Define (with $p\equiv\p{p_\nu^\mu}_{\nu,\mu}$)
\begin{align*}
        \ola\F_{p,r} &= \set{ f\in\ola{E}^\N \Mid
        \forall\mu,\nu\in\N : \ultra f_{p_\nu^\mu,\,r} < \infty } ~,
\\
        \ola\K_{p,r} &= \set{ f\in\ola{E}^\N \Mid
        \forall\mu,\nu\in\N : \ultra f_{p_\nu^\mu,\,r} = 0 }
\\
\text{(resp.}~~
        \ora\F_{p,r} &= \bigcap_{\mu\in\N} \ora\F_{p,r}^\mu
~,~~    \ora\F_{p,r}^\mu = \bigcup_{\nu\in\N}
        \set{ f\in\p{ E_\nu^\mu }^\N \Mid \ultra f_{p_\nu^\mu,r} < \infty }
~,\\
        \ora\K_{p,r} &= \bigcap_{\mu\in\N} \ora\K_{p,r}^\mu
~,~~    \ora\K_{p,r}^\mu = \bigcup_{\nu\in\N}
        \set{ f\in\p{ E_\nu^\mu }^\N \Mid \ultra f_{p_\nu^\mu,r} = 0 }
~\text)~.
\end{align*}
%
%Then, we have the following
%
\begin{propdef}~\parskip-1ex
\begin{enumerate}
\item[(i)] Writing $\olra\cdot$ for both, $\ola\cdot$ or
$\ora\cdot$, we have that $\olra\F_{p,r}$ is an algebra and
$\olra\K_{p,r}$ is an ideal of $\olra\F_{p,r}$; thus, $\olra\G_{p,r}=
\olra\F_{p,r}/\olra\K_{p,r}$ is an algebra.

\item[(ii)] For every $\mu,\nu\in\N,~ d_{p_\nu^\mu}: {(E_\nu^\mu)}^\N
\times {(E_\nu^\mu)}^\N \to\overline\R_+$ defined by
$d_{p_\nu^\mu}(f,g)= \ultra{f-g}_{p_\nu^\mu,r}$ is an
ultrapseudometric on ${(E_\nu^\mu)}^\N$.  Moreover,
$(d_{p_\nu^\mu})_{\mu,\nu}$ induces a topological algebra\footnote
{\label{a} over $(\C^\N,\ultra\cdot_{|\cdot|})$, not over $\C$:
scalar multiplication is not continuous.}
structure on $\ola\F_{p,r}$
% (since $d_{p_\nu^\mu}(0,f\cdot g)\leq
% d_{p_\nu^\mu}(0,f)\,d_{p_\nu^\mu}(0,g)$)
such that the intersection of
the neighborhoods of zero equals $\ola\K_{p,r}$.

\item[(iii)] From (ii), $\ola\G_{p,r}= \ola\F_{p,r}/\ola\K_{p,r}$
becomes a topological algebra (over generalized numbers
$\C_r=\G_{|\cdot|,r}$) whose topology can be defined by the family of
ultrametrics $(\tilde{d}_{p_\nu^\mu})_{\mu,\nu}$ where
$\tilde{d}_{p_\nu^\mu} ([f],[g])=d_{p_\nu^\mu}(f,g)$, $[f]$ standing
for the class of $f$.

\item[(iv)] If $\tau_\mu$ denote the inductive limit topology on
$\F_{p,r}^\mu=\bigcup_{\nu\in\N}((\tilde{E_\nu^\mu})^\N,d_{\mu,\nu})$,
$\mu\in\N$, then $\ora\F_{p,r}$ is a topological algebra\footnotemark[1]
for the
projective limit topology of the family $(\F_{p,r}^\mu,\tau_\mu)_\mu$.

\end{enumerate}
\end{propdef}

%********************************************

We have proved in \cite{DHPV1}:
  \begin{proposition}
\begin{enumerate}
\item[(i)] $\ola\F_{p,r}$ is complete.
\item[(ii)] If for all $\mu\in\N$, a subset of $\ora\F^\mu_{p,r}$ is
bounded iff it is a bounded subset of $\left(E_\nu^\mu\right)^\N$ for
some $\nu\in\N$, then $\ora\F_{p,r}$ is sequentially complete.
\end{enumerate}
\end{proposition}

We showed that various definitions of Colombeau algebras $\bar{\C}$ and $\G$
correspond to the sequence $r_n=1/\log n$, $n\in\N$. The embedding of
Schwartz distributions and of smooth functions into $\G$ is well-known. Also
it is well-known that the multiplication of smooth function embedded into $\G$
is the usual multiplication.
In \cite{DHPV2} we have constructed sequence spaces forming algebras which
correspond to Colombeau type algebras of ultradistributions and periodic
hyperfunctions. The main objective of \cite{DHPV2} was to realise the
embeddings of ultradistribution spaces and periodic hyperfunction spaces into
such algeberas and realise the multiplication of regular elements embedded
into corresponding sequence spaces.

\medskip

In the definition of our sequence spaces $\ora\F_{p,r}$ (resp.
$\ola\F_{p,r}$), we assumed that $r_n$ tends to $0$ as $n$ tends to $\infty$.
%
%(Later, we will have {\bf families} of sequences decreasing to 0.)\\
%
One could consider more general sequences of weights.  But, for example, if
$r_n$ is contained in some compact subset of $(0,+\infty)$ then $\olra E$ can
be embedded in the set-theoretical sense via the canonical map $f\mapsto(f)_n$
($f_n=f$).

If $r_n\to\infty$, $\olra E$ is no more included in $\olra\F_{p,r}$.

%%%% Do not make change in big version %%%
%%%% except ``a discrete topology should be replace by ``the'' .... %%%%

In order to have an appropriate topological algebra containing ``$\delta$'',
we need ``divergent'' sequences; this justifies the choice of $r_n\to0$.
Then, our generalized topological algebra induces the discrete topology on the
original algebra $\olra{E}$.
\\
In some sense, it is an analogy to Schwartz' impossibility statement for
multiplication of distributions~\cite{schw}.

%%%%%%%%%%%%%%%%%%%%%%     SEQUENCES OF SCALES     %%%%%%%%%%%%%%%%%%%%%

\section{Sequences of scales and asymptotic algebras}\label{sect:scales}

\subsection{Sequences of scales}\label{sequences}

We can consider a sequence $\p{r^m}_m$ of positive sequences
$\p{r_n^m}_n$ such that
\[
  \forall m,n\in\N :~ r_{n+1}^m\leq r_n^m ~;~~
  \lim_{n\to\infty}r_n^{m}=0 ~.
\]
In addition to this, we request either of the following conditions :
\begin{eqnarray}
        \forall m,n\in\N :&& r_n^{m+1} \le r_n^{m} \label{II}
\\
\text{or ~ }
        \forall m,n\in\N :&& r_n^{m+1} \ge r_n^{m} ~. \label{I}
\end{eqnarray}

Then let, in the first (resp. second) case:
\begin{eqnarray*}
  \olra\F_{p,r} &=& \bigcup_{m\in\N} \olra\F_{p,r^m} ~,~~
  \olra\K_{p,r} = \bigcap_{m\in\N} \olra\K_{p,r^m} ~
\\
\text{( resp. ~}
  \olra\F_{p,r} &=& \bigcap_{m\in\N} \olra\F_{p,r^m} ~,~~
  \olra\K_{p,r} = \bigcup_{m\in\N} \olra\K_{p,r^m} ~\text) ~,
\end{eqnarray*}
where $p=\p{p_\nu^\mu}_{\nu,\mu}$.

\begin{proposition}
With the previous notations, $\olra\G_{p,r}=\olra\F_{p,r}/\olra\K_{p,r}$ is
an algebra.
\end{proposition}

\begin{proof}
Let us start with the first case (\ref{II}).  $r^{m+1}\le r^m \impl \ultra
f_{r^{m+1}}\ge\ultra f_{r^m}$ if $p(f_n)<1$, hence $\K_{m+1}\subset\K_m$.
Conversely, $\F_{m+1}\supset\F_m$.  Thus, intersection for $\K$
and union for $\F$ makes sense.
Moreover, because of this inclusion property, $\F$ is indeed a subalgebra.
To prove that $\K$ is an ideal, take $(k,f)\in\K\times\F$, i.e.
$\forall m'':k\in\K_{m''}$, and $\exists m':f\in\F_{m'}$.
We have to show that $\forall m: k\cdot f\in\K_m$. So let $m$ be given.\\
If $m<m'$, then $\K_{m'}\subset\K_{m}$, thus
$k\cdot f\in\K_{m'}·\F_{m'} \subset \K_{m'} \subset \K_{m}$.\\
If $m'<m$, then $\F_{m'}\subset\F_{m}$, thus
$k\cdot f\in\K_{m}·\F_{m'} \subset \K_{m}·\F_{m} \subset \K_{m}$.
\\[1ex]%
Now turn to the  the second case (\ref{I}).
The same reasoning gives now $\K_{m+1}\supset\K_m$
and $\F_{m+1}\subset\F_m$, justifying definitions of $\F$ and $\K$.
$\F$ is obviously a subalgebra.
To see that $\K$ is an ideal, take $(k,f)\in\K\times\F$.  Then
$\exists m:k\in\K_m$, but also $f\in\F_m$, in which $\K_m$ is an
ideal.  Thus, $k\cdot f\in\K_m\subset\K$.
%\\[1ex]
\end{proof}

\begin{example}
$
  r_n^{m}=\begin{cases} 1 &\text{if ~} n \leq m\\
                 0 & \text{if ~} n > m \end{cases}
$
(with the convention that $0^0=0$)
gives Egorov--type algebras, where the ``subalgebra'' contains
everything and the ideal contains only stationary null sequences.
\end{example}

% \begin{example}
% $r_n^{m}=1/n^{\frac{m}{m-1}}$ : related to ultradistributions
% \end{example}

\begin{example}
$r_n^m=1/|\log a_m(n)|$, where ${(a_m:\N\to\R_+)}_{m\in\Z}$ is
an asymptotic scale, i.e. $\forall m\in\Z: a_{m+1}=o(a_m)$,
$a_{-m}=1/a_m$, $\exists M: a_M=o(a_m^2)$. This gives back the
asymptotic algebras of~\cite{asymp}, cf. Section~\ref{sect:asympt}.
\end{example}

\subsection{Asymptotic algebras}\label{sect:asympt}

%\subsection{\CEP--algebras}\label{sect:CEP}

Let us recall that \CEP--algebras \cite{CEP} are based on a vector space $\E$
with a filtering family $\P$ of seminorms, and a ring of generalized numbers
$\CC=A/I$. Here, $I$ is an ideal of $A$, which is a subring of $\KK^\Lambda$,
where $\KK=\R$ or $\C$, and $\Lambda$ is some indexing set. Both $A$ and $I$
must be solid as a ring, i.e. $\forall s\in\KK^\Lambda:(\exists r\in A:
\forall\l\in\Lambda: |s_\l| \leq |r_\l|) \Longrightarrow s\in A$, and idem for
$I$. Then, the \CEP--algebra is defined as $ \G_{\CC,\E,\P} = \E_A / \E_I $,
with
\[
  \E_X = \set{ f\in\E^\Lambda \mid \forall p\in\P: p\circ f\in X }
\]
(where $ p\circ f \equiv\lr(){\l\mapsto p(f_\l)}=
\lr(){p(f_\l)}_\l\in(\R_+)^\Lambda \subset \KK^\Lambda$):
In other words, the function spaces $\E_A$ and $\E_I$ are determined by
$\CC=A/I$, by selecting the {\em functions} with the same respective growth
properties than the ``{\em constants\/}''.

It is clear that this is too general to be written in the previously
presented setting of sequence spaces, mainly because there are almost no
%relation between $A$ and
restrictions on $I$.% (one could even have $I=A$).

So let us consider the interesting subclass
%turn to the case
of asymptotic algebras~\cite{asymp}.
Here, $A$ and $I$ are defined by an asymptotic scale\footnote
{
The set $\Lambda$ is supposed to have a base of filters $\B$, to which
the $o(\cdot)$ notation refers to.\\
In the preceding paragraph, %Section~\ref{sect:CEP},
$\forall\l\in\Lambda$ could also be replaced by
$\exists\Lambda_0\in\B,~\forall\l\in\Lambda_0$.
}
$ \mathbf{a} = \left( a_m: \Lambda\to\R_+ \right) _{m\in\Z} $\,:
\begin{eqnarray*}
  A_{\mathbf{a}} &=&\set{ s\in\KK^\Lambda\mid \exists m\in\Z: s=o(a_m) }
\\
  I_{\mathbf{a}} &=&\set{ s\in\KK^\Lambda\mid \forall m\in\Z: s=o(a_m) }
\end{eqnarray*}
Recall that $\mathbf{a}$ must verify: $\forall m\in\N: a_{m+1}=o(a_m)$,
$a_{-m}=1/a_m$, $\exists M: a_M=o(a_m^2)$.
Some examples that have proved to be useful are: %$\Lambda=\N $)

\begin{enumerate}
\item $\Lambda=\N$ and $a_m(\l)=1/\l^m$ :
This leads to Colombeau's generalized numbers and algebras.

\item $\Lambda=\N$ and $a_m(\l)=1/\exp^m(\l)$ for
$m\in\N^*$, where $\exp^m$ is the $m$-fold iterated $\exp$ function:
This gives the so-called exponential algebras~\cite{asymp}.

\item $r_n^{m}=1/n^{\frac{m}{m-1}}$: This is related to
ultradistribution spaces, and is discussed in~\cite{DHPV2}.

\end{enumerate}

\begin{proposition}
Asymptotic algebras can be recovered in our formulation by choosing the
sequence of weights $r^m=1/|\log a_m|$ (i.e. $r^m_\l=1/|\log a_m(\l)|$).
\end{proposition}

\begin{proof}
We will show that $\E_I=\K_{\P,r}$ and $\E_A=\F_{\P,r}$, for
$r^m=1/|\log a_m|$.\\
In view of the definitions, this amounts to show the equivalences
\[
  \forall p, \underset{\displaystyle(\exists)}{\forall} a_m :
   p \circ f = o(a_m)
\iff
  \forall p, \underset{\displaystyle(\exists)}{\forall} r^m :
   \ultra f_{p,r^m} \underset{\displaystyle(<\infty)}{= 0} ~.
\]
%\begin{itemize}
\def\item{}
Let us start with
\item $\E_A\subset\F_{\P,r}:$ Let $f\in\E_A$. Thus,
$\forall p\in\P$, $\exists m: p\circ f=o(a_m)$.
We can assume $a_m > 1$, such that
$ r^m = 1/\log a_m \iff a_m = e^{1/r^m} $.
Thus $ p\circ f = o(e^{1/r^m}) $. But
$ p\circ f < e^{1/r^m} \impl (p\circ f)^{r^m} < e $,
thus $ \limsup (p\circ f)^{r^m} < \infty $
and $f\in\F_{\P,r}$.
%\hfill$\Box$

Conversely,
\item
$ \F_{\P,r} \subset \E_A :$ We have $\forall p\in\P$,
$ \exists \bar m: \limsup (p\circ f)^{1/|\log a_{\bar m}|} < \infty $.
With
$$
  ( p\circ f )^{1/|\log a_m |} \leq C \iff
  p\circ f \le \left( a_m \right) ^{\log C} ~, \quad (a_m,C>1)
$$
we have: $ \exists C>0 $,
$ \exists\Lambda_0:\forall \l\in\Lambda_0:
  p( f_\l ) \le ( a_{\bar m}(\l) )^{|\log C|}
$.
Thus, using the 3rd property of scales,
$\exists m: p\circ f = o(a_m)$.
%\hfill$\Box$

Now turn to
\item $ \E_I \subset \K_{\P,r}:$
We have $ \forall \bar m:p\circ f=o(a_{\bar m}) $.\\
Take $ m \in \N $. Now, for any $ q \in \N $,
$ \exists \hat m : a_{\hat m}=o({a_m}^q) $.
and  $ p\circ f = o(a_{\hat m}) $.

Using $a_m = e^{-1/r^m}$,
$ a_{\hat m} = o( {a_m}^q ) = o((e^{-1/r_m})^q) = o((e^{-q})^{1/r^m}) $,
i.e., $(p\circ f)^{r^m}\le e^{-q}$ for $\l$ ``large enough''.
As $q$ was arbitrary, we have $(p\circ f)^{r^m}\to 0$
and thus $f\in\K$.
%\hfill$\Box$

Finally,
\item $ \K_{\P,r} \subset \E_I :$ We have
$ \forall \bar m: \limsup p(f_\l)^{1/|\log a_{\bar m}|}=0 $,
i.e.,
$$
  \forall C>0,\exists \Lambda_0,\forall \l\in\Lambda_0:
  p(f_\l)^{1/|\log a_{\bar m}|}<C ~.
$$
With $a_m,C<1$, this gives
$p(f_\l)\le C^{|\log a_{\bar m}|}={a_{\bar m}}^{|\log C|}$.
Now, to show that $f\in\E_I$, take any $m$. Let $\bar m=m+1$ and $C=1/e$:
$\exists \Lambda_0, \forall \l\in\Lambda_0 :
p(f_\l)< a_{\bar m}(\l)$.
But $a_{\bar m} = a_{m+1} = o(a_m)$, thus $p\circ f=o(a_m)$.
%\end{itemize}
%\vskip-2ex
\end{proof}

\subsection{Algebras with infra-exponential growth}

A second interesting subclass of \CEP--algebras are of the form
\begin{eqnarray*}
  A &=& \set{ s\in\KK^\Lambda \mid \forall\sigma<0:s=o(a_\sigma) }
\\
  I &=& \set{ s\in\KK^\Lambda \mid \exists\,\sigma>0: s=o(a_\sigma) }
\end{eqnarray*}
where $\mathbf{a}=\left( a_{\sigma}\right) _{\sigma\in\R }$ is again
a scale (i.e. $\forall\sigma>\rho,~a_{\sigma}=o(a_{\rho})$,
etc.), but indexed by a real number.
(Note that here $A$ is given as intersection and $I$ as union of sets,
that's why this case is not covered by the previous one.)\@gobble
{
for this reason I'm not sure whether we can describe this in our
formulation!
}

For example (again with $\Lambda=\N $),
\[
  a_\sigma:=\l\mapsto1/\exp\left( \sigma\,\l\right)
\]
gives the so-called algebras with infra--exponential
growth~\cite{infexp}, pertaining to the embedding of periodic
hyperfunctions in \CEP--algebras.

These algebras can be obtained by taking $\F=\set{f\mid\ultra
f_r\le1}$ and $\K=\set{f\mid\ultra f_r<1}$, with $r_n=\frac1n$. (As
the norm is compared to 1, all scales $r_\sigma=1/|\log a_\sigma|$
(i.e. $r_\sigma(\l)=1/|\sigma\l|$ are equivalent.
More details on this ``dual'' construction, where the conditions
$<\infty$ and $=0$ are replaced by $\le1$ and $<1$, will be
given in a forthcoming %are left to a separate
publication.)

        %%%%%%%%%%%  FUNCTORIAL  %%%%%%%%%%%%

\section{Functorial properties}\label{sect:functor}

In this section, we want to investigate on conditions sufficient to extend
mappings on the topological factor algebras constructed as before. Consider
for example
\[
  \varphi: E\to F
\]
where $(E,P)$ and $(F,Q)$ are spaces equipped
with families of seminorms $P$ and $Q$. 
We shall note in this section $\F_{\Pi,r}(·)$, $\K_{\Pi,r}(·)$ and
$\G_{\Pi,r}(·)$ the spaces defined as above, where · stands for $E$ or
$F$ and $\Pi$ stands for $P$ or $Q$.

Suppose that $\varphi$ satisfies the following hypotheses:
\begin{align*}
  (F_1):&& f\in\F_{P,r}(E) &\impl \vp(f)\in\F_{Q,r}(F)
\\
  (F_2):&& f\in\F_{P,r}(E) ,~ h\in\K_{P,r}(E) &\impl
                        \vp(f+h)-\vp(f)\in\K_{Q,r}(F)
\end{align*}
where we write
$
  \vp(f) := \left( \vp\left( f_n\right) \right) _n
$.
%
%It is clear that this is the general and natural condition for $\vp$
%being well defined on the factor spaces.
%%$\G_{\Pi,r}(·)=\F_{\Pi,r}(·)/\K_{\Pi,r}(·) $
%More precisely,
Then we can consider the following
\begin{definition}
Under the above hypothesis, we define the $r$--extension of $\vp$ by
\[
  \Phi:=\G_r(\vp) :=\left(
  \begin{array}[c]{rcl}
    \G_{P,r}(E) &\to& \G_{Q,r}(F)\\{}
    [f] &\mapsto& \vp(f) +\K_{Q,r}
  \end{array}
\right)
\]
where $f$ is any representative of $[f]=f+\K_{P,r}(E)$.
\end{definition}

The above consideration is of course a very general condition for a
map to be well defined on a factor space. In fact, it does not depend
on details of how the spaces $\F_{P,r}(E)$ and $\K_{P,r}(E)$ are defined.
In particular, here $r$ can also be a family of sequences $\p{r^m}_m$,
and $E$ can be of proj-proj or ind-proj type.

\begin{example}
Consider a linear mapping $u\in\L(E,F)$, continuous for $(P,Q)$. Fix
$q\in Q$. As $u$ is continuous, there exists $p=p_{(q)}$ such that
$$
  \exists c:\forall x\in E :~
  q\left( u(x) \right) \leq c\,p_{(q)}(x) ~.
$$
Thus, $\forall f,h\in E^\N:$
\begin{eqnarray*}
  \limsup \left( p_{(q)} \left( f_n\right) \right) ^{r_n} < \infty
& \impl&
  \limsup\left( q\left( u\left( f_n\right) \right) \right) ^{r_n} < \infty
\\
  \limsup\left( p_{(q)}\left( h_n\right) \right) ^{r_n} = 0
& \impl&
  \limsup\left( q\left( u\left( h_n\right) \right) \right) ^{r_n} = 0
\end{eqnarray*}
\end{example}

This example shows how we can define moderate or compatible maps with
respect to the ``scale'' $r$.  In fact, the concrete definitions will
depend on the monotony properties of the family $(r^m)$ of sequences
of weights, according to which
$\F=\bigcup\F_m$ and $\K=\bigcap\K_m$ (for $r^{m+1}\le r^m$),
or
$\F=\bigcap\F_m$ and $\K=\bigcup\K_m$ (for $r^{m+1}\ge r^m$).
%
%For example, recall that asymptotic algebras correspond to the first case:
%$a_{m+1}=o(a_m) \leadsto \log a_{m+1}<\log a_m \leadsto
%|\log a_{m+1}|>|\log a_m|$ i.e. $r^{m+1}<r^m$ for $r^m=\frac1{|\log a_m|}$.)

The analysis of continuity (in the sense of $\ultra{\cdot}_{p,r}$)
shows that the following definitions are convenient:
\begin{definition}[for $r^{m+1}\le r^m$]
The map $g:\R_+\to\R_+$ is said to be \textbf{$r$--moderate} iff it is
increasing and
\[
  \forall m\in\N ~ \exists M\in\N ~ \forall x\in\R_+ :
  \sup_{n\in\N}\left(g\left(x^{1/r^m_n}\right)\right)^{r^M_n}<\infty
\]
The map $h:\R_+\to\R_+$ is said to be \textbf{$r$--compatible} iff it
is increasing and
\[
  \forall M\in\N ~ \exists m\in\N :
  \left( h \left( x^{1/r^m_n} \right) \right) ^{r^M_n}
  \underset{x\to0}{\longrightarrow}0 \text{~ uniformly in }\,n ~.
\]
\end{definition}

\begin{proposition}
The above definition of an $r$--moderate map $g$ is equivalent to
$$
        g \text{ increasing, and ~} \forall m~\exists M:
        g(\F_{r^m}^+)\subset\F_{r^M}^+  ~.
$$
where $\F^+_{r^m} = \R_+^\N\cap\F_{|·|,r^m}$ are ``moderate'' sequences
of nonnegative numbers.
\\
The definition of an $r$--compatible map $h$ can be written as
$$
        h \text{ increasing, and ~} 
        \forall M~\exists m:\ultra{h(C)}_M\to0\text{ when }\ultra C_m\to0 ~.
$$
or, equivalently
$$
        h\text{ continuous at 0, increasing, and ~}
        \forall M~\exists m:h(\K^+_{r^m})\subset\K^+_{r^M} ~.
$$
\end{proposition}
\begin{proof}
(writing $\F_m$ for $\F^+_{r^m}$):
We have $ g(\F_m)\subset\F_M $
$$
\begin{array}{*9{cr@{~}l}}
\iff&\forall C\in\R_+^\N: C \in\F_m &\impl g(C)\in\F_M \\[.5ex]
\iff&\forall C\in\R_+^\N: \ultra C_m<\infty &\impl\ultra{g(C)}_M<\infty\\[.5ex]
%$. This is equivalent to $
\iff& (\exists x>0,\forall n: C_n<x^{1/r^m_n})
			&\impl \sup_n g(C_n)^{r^M_n}<\infty	\\[.5ex]
\iff&
%$ and finally to $
%&\llap{$\displaystyle
\forall x\in\R_+: \sup_n g(x^{1/r_n^m})^{r^M_n}\kern-1ex&\kern1ex
~<~ \infty ~.
\end{array}$$
%\\
%
For $h$, again take $C_n=x^{1/r^m_n}$, such that
$x\to0\iff \ultra C_m\to0$.\\
Clearly, the first form implies that $h$ is continuous at 0, thus
``$\to0$''
can be replaced by
%equalities.
``$=0$''.
Thus we have
$\forall M~\exists m:C\in\K_m\impl h(C)\in\K_M$, i.e. $h(\K_m)\subset\K_M$.
\end{proof}

\begin{definition}[for $r^{m+1}\ge r^m$]
The map $g:\R_+\to\R_+$ is said to be \textbf{$r$--moderate} iff it is
increasing and
\[
  \forall M\in\N ~ \exists m\in\N ~ \forall x\in\R_+ :
  \sup_{n\in\N}\left(g\left(x^{1/r^m_n}\right)\right)^{r^M_n}<\infty
\]
The map $h:\R_+\to\R_+$ is said to be \textbf{$r$--compatible} iff it
is increasing and
\[
  \forall m\in\N ~ \exists M\in\N :
  \left( h \left( x^{1/r^m_n} \right) \right) ^{r^M_n}
  \underset{x\to0}{\longrightarrow}0 \text{~ uniformly in }n ~.
\]
\end{definition}
\begin{proposition}
The condition of $r$--moderateness can be written
$$
        g\text{ increasing and }
        \forall M,\exists m:g(\F_{r^m}^+)\subset\F_{r^M}^+~.
$$
%(once again only ``almost'' equivalent to $g(\F_r^+)\subset\F_r^+$).
\\
The condition of $r$--compatibility can be written
$$
        h\text{ increasing and }\forall m,\exists M:\ultra{h(C)}_M\to0
        \text{ for }\ultra C_m\to0 
$$
or equivalently
$$
        h\text{ continuous at 0, increasing, and }
        \forall m,\exists M: h(\K_{r^m}^+)\subset\K_{r^M}^+ ~.
$$
\end{proposition}
\begin{proof}
%The proof goes along the same lines than the
The previous proof applies, {\em mutatis mutandis}:
supremum is to be replaced by the corresponding ultranorm,
%replace $\sup(...)^{r^M}$ by $\ultra·_M$ and
%$\forall x,x^{1/r^m}$ by $\forall C\in\F_m$, to get the result.\\
%For $h$, also replace ``uniformly $\to0$'' by $\ultra·_M\to0$
%and $x^{1/r^m},x\to0$ by $C,\ultra C_m\to0$.
and uniform convergence by convergence of the ultranorm.
The definition implies continuity at $h(0)=0$,
so ``$\to0$'' can be replaced by ``$=0$'', 
and thus $h\in\K_M$ ({resp. }$\K_m$).
\end{proof}

\begin{lemma}
If $g$ is $r$--moderate, then $g(\F_r^+)\subset\F_r^+$;
if $h$ is $r$--compatible, then $h(\K_r^+)\subset\K_r^+$.
\end{lemma}
\begin{proof}
Start with the first case, $\F=\cup\F_m$ and $\K=\cap\K_m$:
We have $\forall m~\exists M:g(\F_m)\subset\F_M
\impl \forall m:g(\F_m)\subset\bigcup_M\F_M=\F
\iff \bigcup_m g(\F_m)=g(\F)\subset\F$.\\
Similarly, $\forall M~\exists m:h(\K_m)\subset\K_M
\impl \forall M:\bigcap_m h(\K_m)=h(\K)\subset\K_M
\iff h(\K)\subset\bigcap_M\K_M=\K$.\\
Turn to the second case, $\F=\cap\F_m$ and $\K=\cup\K_m$: The proofs for
$g(\F)$, $h(\K)$ are identical to the proofs for $h(\K)$, $g(\F)$ in the first
case.
\end{proof}

Now we give the definition (valid for both of the above cases)
characterizing maps that extend canonically to $\G_r$:
\begin{definition}
The map $\vp:(E,P)\to(F,Q)$ is said to be \textbf{continuously
$r$--temperate} iff
\begin{align*}
(\alpha)&&
        \exists\,r\text{--moderate } g &,
        \forall q\in Q, \exists p\in P, 
        \forall f\in E: q(\varphi(f)) \le g(p(f))
\\
(\beta)&&
        \exists\,r\text{--moderate } g &,
        \exists\,r\text{--compatible } h:
        \forall q\in Q, \exists p\in P,
\\
&&      \forall f\in E &, \forall k\in E:
        q(\varphi(f+k) - \varphi(f)) \leq g(p(f)) \, h(p(k))
\end{align*}
\end{definition}

\begin{proposition}
Any continuously $r$--temperate map $\varphi$ extends canonically to
\[
        \Phi=\G_r(\vp):\G_{P,r}(E)\to\G_{Q,r}(F) ~.
\]
Furthermore, this canonical extension is continuous for the topologies\\
$\p{\G_{P,r}(E),\p{\smash{\ultra\cdot_{p,r}}}_{p\in P}}$ and
$\p{\G_{Q,r}(F),\p{\smash{\ultra\cdot_{q,r}}}_{q\in Q}}$.
\end{proposition}

\begin{proof}
The proof has two parts: first, the well-definedness of the extension;
secondly, the continuity of $\Phi$.
As a preliminary remark, observe that $\F_{P,r^m}=\set{f\mid\forall p\in P:
p(f)\in\F^+_{r^m}}$, and idem for $\K$. This, and the fact that $\K_{r^m}$
is an ideal in $\F_{r^m}$ (and $\F^+_{r^m}·\K^+_{r^m}\subset\K^+_{r^m}$)
helps us to write the proof using the preceding two characterizations
of moderate and compatible maps.\\
{\em First part of the proof:} We will show that $(\alpha)$ implies
$(F_1)$ and $(\beta)$ gives $(F_2)$.  Using respective definitions of
moderateness and compatibility, the proof will be different for the
two cases $r^{m+1}\le r^m$ and $r^{m+1}\ge r^m$.\\
First case, $r^{m+1}\le r^m$, where $\F=\cup\F_m$ et
$\K=\cap\K_m$:\\
$ad~(F_1):$ Take $f\in\F_{P,r}(E)$, i.e. $\exists m~\forall p:p(f)\in\F^+_m$.
By $(\alpha)$, there is $g$ such that $\exists M:g(\F^+_m)\subset\F^+_M$,
and $\forall q: q(\vp(f))\le g(p(f))\in g(\F^+_m)$, thus 
$\exists M~\forall q:q(\vp(f))\in\F^+_M$, i.e. $\vp(f)\in\F_{Q,r}(F)$.\\[1ex]
$ad~(F_2):$ 
Take $f\in\F$ and $k\in\K$, i.e. $\exists m,\forall p:p(f)\in\F^+_m$
and $\forall m',\forall p:p(k)\in\K^+_{m'}$. Now fix $M$ and $q$.
With $(\beta)$, $\exists g~\forall m~\exists M':g(\F^+_m)\subset\F^+_{M'}$,
and $\exists h~\forall M''~\exists m':h(\K^+_{m'})\subset\K^+_{M''}$.
We use this for $M''=\max(M,M')$, such that $\K^+_{M''}\subset\K^+_{M'}$
and $\K^+_{M''}\subset\K^+_{M}$. Finally, $\exists p: q\p{\vp(f+k)-\vp(f)}
\le g(p(f))\,h(p(k))\in g(\F^+_{m})·h(\K^+_{m'})\subset\F^+_{M'}·\K^+_{M''}$.
Now distingush two cases:
if $M'\le M$, this is in $\F^+_{M'}·\K^+_{M}\subset\F^+_{M}·\K^+_{M}
\subset\K^+_{M}$.
Conversely, if $M<M'$, then this is subset of
$\F^+_{M'}·\K^+_{M'}\subset\K^+_{M'}\subset\K^+_{M}$,
because the $\K^+_{m}$ form a decreasing sequence.
Thus, $\vp(f+k)-\vp(f)\in\K_{Q,r}(F)$.\\[1ex]
Second case, $r^{m+1}\ge r^m$, where $\F=\cap\F_m$ and $\K=\cup\K_m$:\\
$ad~(F_1):$
$f\in\F_{P,r}(E) \iff \forall m~\forall p:p(f)\in\F^+_{m}$.
By $(\alpha)$,
$$\text{
$\exists g~\forall M~\exists m:g(\F^+_m)\subset\F^+_M$,
and $\forall q~\exists p: q(\vp(f))\le g(p(f))\in g(\F^+_m)$.
}$$
Thus, $\forall M~\forall q:q(\vp(f))\in\F^+_M$,
i.e. $\vp(f)\in\F_{Q,r}(F)$.\\[1ex]
$ad~(F_2):$ 
Take $f\in\F$ and $k\in\K$, i.e. $\forall m,\forall p:p(f)\in\F^+_m$
and $\exists m',\forall p:p(k)\in\K^+_{m'}$. Now fix $q$.
With $(\beta)$,
$$\text{
$\exists h~\forall m'~\exists M: h(\K^+_{m'})\subset\K^+_M$
and $\exists g~\forall M~\exists m:g(\F^+_m)\subset\F^+_{M}$,
}$$
and there exists $p$ such that
$$\text{
 $%\exists p:
q\p{\vp(f+k)-\vp(f)} \le g(p(f))\,h(p(k))
\in g(\F^+_{m})\,h(\K^+_{m'})\subset\F^+_{M}·\K^+_{M}\subset\K^+_M$,
}$$
thus $\vp(f+k)-\vp(f)\in\K_{Q,r}(F)$.\\[1ex]
{\em Second part of the proof : continuity of $\Phi$.}
We must show that
$$
        \forall q\in Q:\ultra{\vp(f+k)-\vp(f)}_{q,r^M}\to 0
\text{~ when ~}
        \forall p\in P:\ultra k_{p,r^m}\to 0
$$
%and this $\
for all $M$ (resp. for some $M$), in respective cases.
The proof goes analogous to the above proof of $(F_2)$, by
replacing $p(f)\in\F^+_m$ with $% \leadsto
\ultra f_{p,m}\le K$,
and
$p(k)\in\F^+_m$ with $%\leadsto
\ultra k_{p,m}\le \veps$, etc.
\end{proof}

        %%%%%%%%%%% ASSOCIATION %%%%%%%%%%%%

\section{Association in $\G$}\label{sect:association}

We will introduce different types of association, according to what
has already been considered in the literature on generalized function
spaces.  Generally speaking, we will adopt the following terminology:
{\em strong association\/} is expressed directly on the level of the
factor algebra, while {\em weak association\/} will be defined in
terms of a duality product, and thus with respect to a certain test
function space.
% and using representatives of the algebra's elements.

%\subsubsection*
\paragraph{Association in Colombeau type generalized numbers.}
To start with, recall that Colombeau generalized numbers $[x]$
and $[y]$ are said to be associated, $[x]\approx[y]$, iff
\[
  x_n-y_n \underset{n\to\infty}{\longrightarrow}0
  \text{ ~ ( in $\C$ ) .}
\]
This can also be expressed by considering the subset of null
sequences, $N=\set{x\in\C^\N\mid\lim x_n=0}$, and by defining
% the association by
$[x]\approx [y]\iff x-y\in N$

As any element $j$ of the ideal verifies $j_n\to0$,
this is clearly independant of the representative.
In other words, it is well defined because $I\subset N$.

\subsection{The general concept of $\J,X$--association}

The following general concept of association allows to recover
all known notions of association, as well as the the types we
shall consider below.

\begin{definition}[$\J,X$--association]
Let $\J$ be an additive subgroup of $\F$ containing the ideal $\K$,
and $X$ a set of generalized numbers. Then,
two elements $F,G\in\G=\F/\K$ are called $\J,X$--associated,
$$
        F \underset{\J\!,X}\approx G \text{~ iff ~}
        \forall x\in X: x\cdot(F-G) \in \J/\K ~.
$$
For $X=\set1$, we simply write
$$ F \underset{\J}\approx G \iff F-G \in \J/\K ~.$$ 
\end{definition}

\begin{remark}
As $\J$ is not an ideal, association is not compatible with
multiplication in $\F$ (not even by generalized numbers, only by
elements of $E$).
However, in the case of differential algebras, $\J$ is usually chosen
such that $\smash{\underset{\J,X}\approx}$ is stable under
differentiation.
\end{remark}

\begin{example}
Usual association of generalized numbers, as recalled above,
is obtained for $\J=N$, the set of null sequences:
$$
        [x]\approx[y] \iff [x]\underset{N}\approx[y] ~.
$$
As already mentioned, all elements of the ideal $\K$ tend to zero,
i.e. $\K\subset N$,
as needed for well-definedness at the level of the factor algebra.
\end{example}

\subsection{Strong association}

As mentioned, strong association is defined directly in terms of
the ultranorm (or ultrametric) of elements of the factor space.

\begin{definition}
For $s\in\R_+$, {\bf strong $s$--association} is defined by
$$
        F\overset{s}\simeq G \iff
        \smash{F \underset{\J^{(s)}_{\P,r}}\approx G}
$$ with
\begin{align}
        \J^{(s)}_{\P,r}
        =\set{ f\in\F \mid \forall p\in\P: \ultra f_{p,r}<e^{-s}} ~,
\label{J(s)}
\end{align}
which is equivalent to say
\[
        F\overset{s}\simeq G \iff \widetilde{d}_{p,r}(F,G)<e^{-s} ~.
\]
For $s=0$, we write $F\simeq G$ and simply call them {\bf strongly 
associated}.
\end{definition}

If one has $ F\overset{s}\simeq G$ for all $s\ge0$, then $F=G$.
Indeed, this means that $F-G$ is in the intersection of all balls
of positive radius, which is equal to $\K=0_{\G}$.

\subsection{Weak association in $\protect\olra\G_{p,r}$}

In contrast to the above, weak association is defined by comparing
sequences of {\em numbers\/} (not {\em functions\/}), obtained by
means of a duality product
$$
        \langle\cdot,\cdot\rangle:\olra E\times\DD\to\C ~,
$$
where $\DD$ is a test function space such that $E\hookrightarrow\DD'$
(as for example $\DD=\D$ for $E=\CC^\infty$).
The subset $\J$ defining the association will then be of the form
\begin{align}
        \J = \J_M=\set{ f\in \olra E^\N \mid
                 \forall\psi\in\DD: \p{\langle f_n,\psi\rangle}_n\in M } ~,
\label{J_M}
\end{align}
where $M$ is some additive subgroup of $\C^\N$ (like e.g. $M=N$, the null
sequences).

\begin{example}
For the choices given above, $\DD=\D$, $E=\CC^\infty$ and $M=N$, in
the case of Colombeau's algebra, we get the usual, so-called weak
association $[f]\approx[g]\iff f_n-g_n \to 0$ in $\D'$.

Again, this is independent of the representatives, because $\J\supset\K_{r,p}$.
To see this, consider $j\in\K_{r,p}$. Then for any $\e>0$
there is $n_0$ such that for $n > n_0$,
\[
  |\,\langle j_n,\psi\rangle \,| \leq \e^{1/r_n}\int|\phi|
  \underset{n\to\infty}{\longrightarrow} 0~.
\]
Thus,
$
  \langle f_n,\psi\rangle \underset{n\to\infty}{\longrightarrow}0 \iff
  \langle f_n+j_n,\psi\rangle \underset{n\to\infty}{\longrightarrow}0
$.

This is a special case of the following definition.
\end{example}

\begin{definition}
{\bf$s-\DD'$--association} is defined by
$$
        F\overset s{\underset\DD\approx} G
        \iff \smash{}{F\underset{\J_N,X_s}\approx G}
$$
with $X_s = \set{\left[\p{e^{s/r_n}}_n\right]}$ for $s\in\R$.
\end{definition}
Note that this generalized number is always of the same form, but
depends in each case on the sequence $\p{r_n}_n$ defining the
topology.

\begin{example}
In Colombeau's case, $r=1/\log$, we have $X_s = \set{\left[\p{n^s}_n\right]}$.
For $s=0$ ($X_0=\set1$), we get the already mentioned weak association.
%of Colombeau's algebra.

For $s\ne0$, $[f]\overset s{\underset\D\approx}[g]\iff n^s(f_n-g_n)\to0$
in $\D'$.
This also has already been considered (with $\DD=\D$),
for example in \cite{CEP} (where it had been denoted by $\underset s\approx$).
This association is of course stronger than the simple weak association
(again, because association is not compatible with multiplication
even only by generalized numbers).

As an extension of this example, consider $\J$ as above, and\\
$X=\set{\left[\p{n^s}_n\right]}_{s\in\N}$. This means that
$$
        [f]\approx[g]\iff\forall s\in\N:\lim n^s(f_n-g_n) = 0 \text{ in } \D'
$$
In Colombeau's algebra, this amounts in fact to strict equality.

%The same constructions can be applied to examples~\ref{Ex-Sobolev}
%and~\ref{Ex-Full-Colombeau} (generalized Sobolev space and full
%Colombeau algebra.
%In the case of ultradistributions, we take $\DD=\D^{(m)}$ and
%$e^{s/r_n} = \exp[s\,n^{\frac1{m'-1}}]$ for Beurling case, and
%analogous definitions in the Roumieu case. 
%For periodic hyperfunctions (with $\DD=\A(\T)$) this is also a new
%construction.
\end{example}

%\subsubsection{Weak $s$--association}

\begin{definition}
\textbf{Weak $s$--association} is defined for any $s\in\R$ by
$$
        F{\overset{(s)}\simeq} G \iff \smash{}F\underset{\J_{(s)}}\approx G
$$
where
\begin{eqnarray*}
        \J_{(s)} &=& \set{ f\in E^\N\mid \forall\psi\in\DD:
        \ultra{\p{\langle f_n-g_n,\psi\rangle}_n}_{|\cdot |,r} < e^{-s}}
\\      &=& \set{ f\in E^\N\mid \forall\psi\in\DD:
        \limsup_{n\to\infty}\lr||{\langle f_n-g_n,\psi\rangle}^{r_n}<e^{-s}}
        ~.
\end{eqnarray*}
It is obtained from the general setting (\ref{J_M}) by observing that
$\J_{(s)}=J_M$ with (cf. eq. (\ref{J(s)}))
$$
        M = \J_{|\cdot|,r}^{(s)}
        = \set{ c\in\C^\N \mid \ultra c_{|\cdot |,r}< e^{-s}} ~.
$$
%(Let us note, for the sake of completeness, that this is an ideal in the 
%subalgebra $H = \set{ c\in\C^\N \mid \ultra c_{|\cdot|,r} \le 1 }$.)

For $s=0$, we write $F\overset{\rm sw}\approx G$
and call $F$ and $G$ {\bf strong--weak associated}. 
\end{definition}

%\begin{remark}
Weak $s$--association implies $s-\DD'$--association, but conversely
$s-\DD'$--association only implies weak $s'$--association for all $s'<s$.
%\end{remark}

%\begin{example} ?
%\end{example}

%\begin{remark}
Let us consider some details concerning the structure of strong-weak
association. In the following we will note $|·|_r=\ultra·_{|·|,r}$,
i.e.
\[
        |c|_r = \limsup_{n\to\infty} |c_n|^{r_n} ~.
\]
To start with, let us remark that
$\I_{|\cdot|,r} =\set{ c\in\C^\N \mid | c |_r < 1 }$
is an ideal in the subalgebra
$\H_{|\cdot|,r} =\set{ c\in\C^\N \mid | c |_r \le 1 }$ of $\C^\N$.

Let us now consider the topology on $\C^\N$ induced by the
$|\cdot|_r$--norm. We have
\[
        |c|_r \leq a  \iff
        \forall b>a, \exists n_0, \forall n>n_0 :
        |c_n| \leq b^{1/r_n} ~.
\]
But now observe that for $b>1$, $b=1$ and $b<1$, the limit of the
last expression is respectively $\infty$, $1$ and $0$. This means that
\begin{enumerate}
\item $|c|_r < 1 \impl \lim c_n = 0$

\item $\lim c_n=0 \impl \forall b>1~ \exists n_0~ \forall n\geq n_0:$
  $ |c_n|\leq b^{1/r_n}\to\infty, \impl |c|_r \leq 1$

\item $|c|_r=1:$ any value in $\R_+\cup\set{\infty}$ (or none at all)
is possible as limit for $c_n$. Indeed, whatever be the null sequence
$(r_n)$, the sequences $c_n=r_n$ (resp. $c_n=1/r_n$) have limits $0$ 
(resp. $\infty$), but
\[
  |c_n|^{r_n}= \exp(\pm r_n\log r_n) \underset{n\to\infty}\to 1
~~~
  (\text{ because } x\log x\underset{x\to0}\to 0 ~)
\]
i.e. $|c|_r=1$.
\end{enumerate}

Thus, all elements of the open unit ball are weakly associated to
zero. This is very similar to classical results related to ultrametric
spaces and weak topology.
%\end{remark}

\begin{proposition}
Weak $s$--association implies $s-\DD'$--association,
but the converse is not true.
\end{proposition}
\begin{proof}
This follows from $ \ultra c_{|\cdot|,r} < 1 \impl
\lim\limits_{n\to\infty} c_n = 0 \impl
\ultra c_{|\cdot|,r} \leq 1 $, with $c_n=\langle
f_n,\psi\rangle\,e^{s/r_n}$.
%\\
%
(As already seen, for $\ultra c_{|\cdot|,r} = 1$,
nothing can be concluded about
%eveything is possible for
the limit of $(c_n)$).
\end{proof}

\bigskip\parindent0pt\parskip1ex\flushleft

Antoine Delcroix: IUFM de la Guadeloupe et Laboratoire A.O.C.,\\
Morne Ferret, BP 399, 97159 Pointe \`a Pitre cedex (Guadeloupe, F.W.I.)\\
Tel.: 00590 590 21 36 21, Fax : 00590 590 82 51 11,\\
E-mail: {\tt Antoine.Delcroix@univ-ag.fr}

Maximilian Hasler: Universit\'e des Antilles et de la Guyane,\\
Laboratoire A.O.C., D\'epartement Scientifique Interfacultaire,\\
BP 7209,
97275 Schoelcher cedex (Martinique, F.W.I.)\\
Tel.: 00596 596 72 73 55, Fax : 00596 596 72 73 62,\\
e-mail : {\tt Maximilian.Hasler@martinique.univ-ag.fr}

Stevan Pilipovi\'c: University of Novi Sad, Institute of Mathematics, \\
Trg D. Obradovi\'ca 4, 21000 Novi Sad (Yougoslavia)\\
Tel.: 00381 21 58 136, Fax: 00381 21 350 458,\\
e-mail: {\tt pilipovic@im.ns.ac.yu}

Vincent Valmorin: Universit\'e des Antilles et de la Guyane,\\
Laboratoire A.O.C.,  D\'epartement de Math\'ematiques, \\
Campus de Fouillole, 
97159 Pointe \`a Pitre cedex (Guadeloupe, F.W.I.)\\
Tel.: 00590 590 93 86 96, Fax: 00590 590 93 86 98,\\
e-mail: {\tt Vincent.Valmorin@univ-ag.fr}

\end{document}